\documentclass[numreferences]{kluwer}
\usepackage{amssymb,latexsym,eucal,mathrsfs}

\usepackage[curve,v2]{xypic}

\newcommand{\tensor}{\otimes}

\newcommand{\cliff}{\bf{Cliff}}
\newcommand{\picard}{\bf{Pic}}

\newcommand{\wt}[1]{\widetilde{#1}}
\newcommand{\bt}{\boxtimes}
\newcommand{\wtp}{\widetilde{\P^n}}

\newcommand{\blow}[2]{{\rm Bl}_{#2}{(#1})}

\newcommand{\ses}[3]{0\rightarrow#1\rightarrow#2
   \rightarrow#3\rightarrow0}

\newcommand{\K}{{\mathscr K}}

\newcommand{\F}{{\mathcal F}}

\newcommand{\E}{{\mathscr E}}

\renewcommand{\I}{{\mathscr I}}

\renewcommand{\O}{{\mathcal O}}
\renewcommand{\P}{{\mathbb{P}}}
\newcommand{\C}{{\mathbb{C}}}

\newcommand{\Z}{{\mathbb{Z}}}

\newenvironment{proof}{\par \medskip \noindent
{\sc Proof:}}{\nopagebreak \hfill $\Box$ \par \medskip}

\begingroup
\newtheorem{thm}{Theorem}[section]   
\newtheorem{cor}[thm]{Corollary}     
\newtheorem{lemma}[thm]{Lemma}         
\newtheorem{prop}[thm]{Proposition}  
\newtheorem{conj}[thm]{Conjecture}        
   
\endgroup

\newenvironment{rem}[2]{\refstepcounter{thm} \label{#2} 
\par \medskip \noindent {\bf #1 \thethm }}{\par \medskip}

\begin{document}

\begin{article}
\begin{opening} 

\author{Peter \surname{Vermeire}\email{petvermi@math.okstate.edu}} 
\runningauthor{Peter Vermeire}
\institute{Oklahoma State University}
\begin{ao}
Department of Mathematics,
Oklahoma State University, 
Stillwater OK 74078
\end{ao}

\title{On the Regularity of Powers of Ideal Sheaves}





\begin{abstract}
We use the geometry of the secant variety to an embedded smooth curve
to prove some vanishing and regularity theorems for powers of ideal
sheaves.
\end{abstract}

\classification{1991 Subject Classification}{Primary
14F17; Secondary 14H60,14N05}

\end{opening}


\section{Introduction}
The original motivation for this work was an attempt to understand
an unpublished manuscript of J. Rathmann in which he
proves the following non-trivial result via a fairly lengthy
calculation on the triple product $C\times C\times C$:
\begin{thm}\cite{rathmann}\label{rathorig}
Let $C\subset\P^n$ be a smooth curve embedded by a line bundle of
degree at least $2g+3$.  Then $H^i(\P^n,\I^2_C(k))=0$ for $k\geq 3$,
$i>0$. 
{\nopagebreak \hfill $\Box$ \par \medskip}
\end{thm}

This result was used by A. Bertram to obtain the following: 
\begin{thm}\label{bertramorig}\cite[4.2]{bertramvanishing}
Let $C\subset\P^n$ be a smooth curve embedded by a line bundle of
degree at least $\frac{8g+2}{3}$.  Then $H^i(\P^n,\I_C^a(k))=0$ for
$k\geq 2a-1$, $i>0$.
{\nopagebreak \hfill $\Box$ \par \medskip}
\end{thm}
Bertram proceeds quite differently than Rathmann: using the GIT flip construction of
Thaddeus \cite{thad}, as well as results from
\cite{bertramhard},\cite{bertramlog}, he constructs useful log canonical
divisors on the blow up of $\P^n$ along $C$, and then obtains
vanishing results from a Kodaira-type vanishing theorem.

We work in the same general context as Bertram, though we mostly avoid
the explicit use of flips and of generalized Kodaira-type vanishing, to give a new proof
(Corollary~\ref{simplevanishing}) of
Rathmann's result and then to prove an extension of
Theorem~\ref{bertramorig} suggested in \cite{bertramvanishing}:
\begin{thm}\label{introthm}
Let $C\subset\P^n$ be a smooth curve embedded by a line bundle of
degree at least $2g+3$.  Then $H^i(\P^n,\I_C^a(k))=0$ for
$k\geq 2a-1$, $i>0$.
{\nopagebreak \hfill $\Box$ \par \medskip}
\end{thm}
See Theorem~\ref{generalvanishing} for a slightly more general statement.

As the title suggests, we use these to make statements regarding the
regularity of powers of ideal sheaves (Corollary~\ref{simplevanishing}). 
We also include some statements for canonical curves
(Proposition~\ref{canonical}) which should not be considered at all optimal.  

Finally, we mention a closely related conjecture of Wahl, toward which
we hope to adapt these techniques:
\begin{conj}\cite{wahl}
Let $C\subset\P^{g-1}$ be a canonically embedded curve with $\cliff C\geq
3$. Then $H^1(\I^2_C(3))=0$.
{\nopagebreak \hfill $\Box$ \par \medskip}
\end{conj}

Note that many of the results in Section 2
(Corollary~\ref{arbdeghyp} though Proposition~\ref{canonical}) can be
derived from more general results on point sets
\cite{chandler2},\cite{ggp}.  Further, these results, along with
Rathmann's Theorem~\ref{rathorig}, can be used to derive
Theorem~\ref{generalvanishing} (and hence
Theorem~\ref{introthm}).  However, we have
chosen to retain these results and their proofs as they give  
context to the main results and illustrate that fact that only the
$k=3$ statement in Theorem~\ref{rathorig} is not elementary.
Specifically, the proofs in Section 2 are quite short and are in much
the same spirit as the proof of the main result, Theorem~\ref{rathext}. 

{\bf Acknowledgements:} I would like to thank Sheldon
Katz, Zhenbo Qin, and Jonathan Wahl for their helpful conversations
and communications.  I would especially like to thank the referee for
pointing out some relevant results on point sets that I was not,
though should have been, aware.

\section{Elementary Vanishing}
In this section, we collect a few fairly elementary vanishing
statements that are somewhat broader than those mentioned in the
introduction.  The first is due to Lazarsfeld (Cf. \cite[2.3]{wahl2}):
\begin{lemma}\label{killconormal}
Let $C\subset \P^n$ be a smooth curve, scheme theoretically defined by
forms of degree $r$.  If $H^1(C,\O_C(t))=0$, then
$H^1(C,N^{*^{\tensor a}}(k))=0$ for $k\geq ra+t$. 
{\nopagebreak \hfill $\Box$ \par \medskip}
\end{lemma}
We do not present a proof, as we will next describe a direct extension
of the technique.  Let $C\subset \P^n$ be a smooth
curve scheme theoretically cut out by hypersurfaces of degree $r$.
Tensor the resolution of the ideal sheaf:
$$\rightarrow \bigoplus_i\O_{\P^n}(-a_i)\rightarrow
\Gamma(\I_C(r))\tensor \O_{\P^n}(-r)\rightarrow \I_C\rightarrow 0$$ 
by $\O_C(k)$, and break the sequence into two short exact sequences:
$$\ses{\K_2}{\bigoplus_i\O_C(k-a_i)}{\K_1}$$
$$\ses{\K_1}{\Gamma(\I_C(r))\tensor \O_C(k-r)}{N^*_C(k)}$$
Suppose $a_{i}\geq a_{i+1}$ for all $i$.
If $H^1(\O_C(k-a_1))=0$, then $\Gamma(\I_C(r))\tensor
\Gamma(\O_C(k-r))\rightarrow \Gamma(N^*_C(k))$ is surjective (note
that $H^1(\O_C(k-a_i))=0$ because $C$ is a curve and by maximality of $a_1$).
Now if $H^1(\I_C(k-r))=0$, we have the diagram
\begin{center}
{\begin{minipage}{1.5in}
\diagram
 \Gamma(\I_C(r))\tensor \Gamma(\O_{\P}(k-r)) \dto \rto &
 \Gamma(\I_C(k)) \dto & \\
 \Gamma(\I_C(r))\tensor \Gamma(\O_C(k-r)) \dto \rto & \Gamma(N^*_C(k))
 \rto & 0\\
 0 & & 
\enddiagram
\end{minipage}}
\end{center}
where the first vertical map is surjective by the normality hypothesis
just mentioned, and the second horizontal map is surjective by the
above discussion.  Therefore, the second vertical map is surjective
and we have:
\begin{cor}\label{arbdeghyp}
Let $C\subset \P^n$ be a smooth curve scheme theoretically cut out by
hypersurfaces of degree $r$, with syzygies generated by forms of
degree at most s.  If $H^1(\I_C(k-r))=H^1(\O_C(k-r-s))=0$,
then $H^1(\I^2_C(k))=0$.
\end{cor}
\begin{proof}
This follows by the above discussion and the sequence 
$$\ses{\I^2_C(k)}{\I_C(k)}{N^*_C(k)}$$
(note that the vanishing of $H^1(\I_C(k))$ is implied by our
assumptions).
\end{proof}

\begin{rem}{Terminology}{staywithnicecurves}
For the remainder of the paper, we will be interested in curves that
are at least projectively normal and whose homogeneous ideals are
generated by quadrics.  This is usually referred to as Green's
condition $(N_1)$.  If, further, the syzygies among the defining
quadrics are generated by linear relations, we have condition $(N_2)$.
Recall that a smooth curve embedded by a line bundle of degree at
least $2g+1+p$ satisfies condition $(N_p)$ \cite{mgreen}.
{\nopagebreak \hfill $\Box$ \par \medskip}
\end{rem}

Proceeding inductively, we just as easily deduce 
vanishing statements for higher powers of the ideal sheaf.  We will
assume, however, that $H^1(\O_C(1))=0$; the more general case may be
similarly worked out.  In particular, tensoring the resolution of the
ideal by $S^aN^*_C(2a+1)$ and applying Lemma~\ref{killconormal}, we
see that the map $$\Gamma(\I_C(2))\tensor
\Gamma(S^aN^*_C(2a+1))\rightarrow \Gamma(S^{a+1}N^*_C(2a+3))$$ is
surjective.  There is an analogous diagram
\begin{center}
{\begin{minipage}{1.5in}
\diagram
 \Gamma(\I_C(2))\tensor \Gamma(\I^a_C(2a+1)) \dto \rto &
 \Gamma(\I^{a+1}_C(2a+3)) \dto & \\
 \Gamma(\I_C(2))\tensor \Gamma(S^aN^*_C(2a+1)) \dto \rto &
 \Gamma(S^{a+1}N^*_C(2a+3)) \rto & 0\\
 0 & & 
\enddiagram
\end{minipage}}
\end{center}
where the first vertical map is surjective by the previous stage.  
\begin{prop}\label{basic}
Let $C\subset \P^n$ be a smooth curve with $H^1(\O_C(1))=0$.
\begin{enumerate}
\item If $C$ satisfies $(N_1)$, then:
\begin{enumerate}
\item $H^1(\I^a_C(k))=0$ for $k\geq 2a+1$
\item $H^2(\I^a_C(k))=0$ for $k\geq 2a-1$ hence $\I^a_C$
is $(2a+2)$-regular 
\end{enumerate}
\item If $C$ satisfies $(N_2)$, then:
\begin{enumerate}
\item $H^1(\I^a_C(k))=0$ for $k\geq 2a$ 
\item $\I^a_C$ is $(2a+1)$-regular
\end{enumerate}
\end{enumerate}
\end{prop}

\begin{proof}
The first part follows directly from the discussion above and
Lemma~\ref{killconormal}.  The second part follows from the fact that the
$\O(-4)$ term in the second stage of the resolution of the ideal may be
removed.
\end{proof}


For the sake of completeness, we include a result not covered by
the above statements, but which may be of interest:
\begin{prop}(Cf. \cite{wahl})\label{canonical}
Let $C\subset\P^{g-1}$ be the canonical embedding of a smooth curve
with $\cliff C\geq 3$.  Then $H^1(\I^a_C(k))=0$ for $k\geq 2a+1$ and
$H^2(\I^a_C(k))=0$ for $k\geq 2a-1$.
\end{prop}

\begin{proof}
This follows exactly as above taking into account:
\begin{enumerate}
\item $H^1(\O_C(2))=0$
\item $C\subset\P^{g-1}$ satisfies condition $(N_2)$
(\cite{schreyer},\cite{voisin}) 
\item $H^1(S^aN^*_C(k))=0$ for $k\geq 2a+1$ (\cite[Thm 2]{belwahl},
\cite{sonny}) 
\end{enumerate}
\end{proof}





We conclude this section by recalling a pair of basic lemmas; the
first describes some situations where the cohomology of powers of
ideal sheaves vanishes ``automatically'', the second gives the
relationship between powers of ideal sheaves and divisors on the
blow-up along the subvariety.

\begin{lemma}\label{elem}
Let $X\subset \P^n$ be a non-degenerate smooth variety of dimension
$r$.  Then
\begin{enumerate}
\item $H^i(\P^n,\I_X^a(k))=0$ for $i\geq r+2$ and $a,k\geq 1$
\item If $H^{r+1}(\P^n,\I_X^a(k))=0$ then
$H^{r+1}(\P^n,\I_X^a(k+\sigma ))=0$ for $\sigma\geq 0$, $a\geq 1$
\item $H^0(\P^n,\I_X^a(k))=0$ for $k\leq a$
\end{enumerate}
\end{lemma}

\begin{proof}
The first two statements follows immediately from the basic sequence
$$\ses{\I_X^{a+1}}{\I_X^a}{S^aN^*_X}$$
The third is just the statement that a form of degree $k$ cannot
vanish $k$ times on a non-degenerate variety.
\end{proof}

\begin{lemma}\label{basicbel}\cite[1.2,1.4]{bel}
Let $Y\subset X$ be a smooth subvariety of codimension $e$ of a smooth
projective variety, $L$ an
invertible sheaf on $X$, $\pi:B=\blow{X}{Y}\rightarrow X$ the blow up
along $Y$ with exceptional
divisor $E$.
Then 
\begin{enumerate}
\item If $0\leq t\leq e-1$, then $H^i(B,\pi^*L(tE))=H^i(X,L)$, $\forall i$
\item $\pi_*\pi^*L(-kE)=\I_Y^k\tensor L$ and $R^i\pi_*(\pi^*L(-kE))=0$
for $k\geq 0$, $i>0$; hence
\item $H^*(X,\I_Y^k\tensor L)=H^*(B,\pi^*L(-kE))$, $k\geq 0$
{\nopagebreak \hfill $\Box$ \par \medskip}
\end{enumerate}
\end{lemma}

\section{The Square of the Ideal Sheaf}

This section is devoted to the proof of Theorem~\ref{rathorig} stated in
the Introduction.
We denote the $i^{th}$ secant variety to an embedded curve $C$ by
$\Sigma_iC$, or just by $\Sigma_i$ when no confusion will result.
The following construction and ``Terracini Recursiveness'' result 
of A. Bertram provides the means for our vanishing results.
Recall that a line bundle $L$ on a curve $C$ is said to be {\em
$k$-very ample} if $h^0(C,L(-Z))=h^0(C,L)-k$ for all $Z\in S^kC$.

Let $C\subset X_0=\P(H^0(C,L))$ be a smooth curve embedded by a
$2k$-very ample line bundle $L$.  Construct a birational morphism
$f:\wt{X} \rightarrow X_0$ which is a composition of the following
blow-ups: 

$f^1:X_1\rightarrow X_0$ is the blow up of $X_0$ along $C=\Sigma_0$

$f^2:X_2 \rightarrow X_1$ is the blow up along the proper transform of
$\Sigma_1$ 

\hskip .5in $\vdots$

$f^k:\widetilde X = X_k \rightarrow X_{k-1}$ is the blow up along the
proper transform of $\Sigma_{k-1}$

We then have:
\begin{thm}\cite[Theorem
1]{bertramhard},\cite[3.6]{bertramlog}\label{terracini} 
Hypotheses as above:
\begin{enumerate}
\item For $i\leq k-1$, the proper transform of each $\Sigma_i$
in $X_i$ is smooth and irreducible of dimension $2i+1$, transverse to
all exceptional divisors, and so in particular $\widetilde X$ is smooth. 
Let $E_i$ be the proper transform in $\widetilde X$ of 
each $f^i$-exceptional divisor.  Then $E_1 + ... + E_k$ is 
a normal crossings divisor on $\wt{X}$ with $k$ smooth components.

\item  
Suppose $i\leq k-1$ and $x\in \Sigma_i\setminus \Sigma_{i-1}$. Then
the fiber $(f^i)^{-1}(x) \subset X_i$ is naturally isomorphic to
$\P(H^0(C,L(-2Z)))$, where $Z$ is the unique divisor of degree $i+1$
whose span contains $x$.  Moreover, the fiber $f^{-1}(x) \subset E_{i+1}
\subset \wt{X}$ is isomorphic to $\widetilde X_Z$, 
the variety obtained by applying the above construction to
the line bundle $L(-2Z)$.
\end{enumerate} 
{\nopagebreak \hfill $\Box$ \par \medskip}
\end{thm}

Write $\picard \wt{X}=\Z H+\Z E_1+\ldots +\Z E_k$. We collect a few 
technical implications of the construction:

\begin{lemma}\label{tech}
Hypotheses and notation as above:
\begin{enumerate}
\item $\Sigma_i\subset X_0$ is normal for $i\leq k-1$
\item $f:E_1\rightarrow C$ is a smooth morphism
\item $f_*\O_{\wt{X}}=\O_{X_0}$ and $R^jf_*\O_{\wt{X}}=0$ for $j\geq 1$
\item $f_*\O_{E_i}=\O_{\Sigma_{i-1}}$
\end{enumerate}
\end{lemma}

\begin{proof}
The first statement follows from the fact that $E_i$ is smooth and
$f:E_i\rightarrow \Sigma_{i-1}$ has reduced, connected fibers.
The second follows from the smoothness of $C$ and the description of
the fibers $f^{-1}(x)\cong \widetilde X_Z \subset E_1$
above. 

For the third, $f_*\O_{\wt{X}}=\O_{X_0}$ is Zariski's Main Theorem.  
$R^jf_*\O_{\wt{X}}=0$ follows from the fact that $H^j(\wt{X},f^*\O_X(mH))=0$
for $m,j>0$ (Cf. \cite[2.69]{kollar-mori}).

To show $f_*\O_{E_i}=\O_{\Sigma_{i-1}}$, note that $f:E_i\rightarrow
\Sigma_{i-1}$ is the composition of birational morphisms to smooth
varieties, followed by a projective bundle, followed by a birational
morphism to $Sec^{i-1}C$ which is normal by the first statement.
\end{proof}

We recover Rathmann's result (Corollary~\ref{simplevanishing}) from the main
result of this section: 

\begin{thm}\label{rathext}
Let $C\subset \P^n$ be a smooth curve embedded by a non-special 
line bundle $L$.  Suppose there exists a point $p\in C$ such that 
$L(2p)$ is $6$-very ample and such that $C\subset \P H^0(C,L(2p-2q))$
satisfies condition $(N_2)$ for all $q\in C$.
Then  $H^i(\P \Gamma(C,L(2p-2q)),\I_C^2(k))=0$, $k\geq 3$, $i>0$.
\end{thm}

The idea of the proof is to take weak (i.e. asymptotic) vanishing
statements on the spaces $X_i$ from Bertram's construction and to
descend them to effective vanishing results along the fibers of
$f:\wt{X}\rightarrow X_0$. For this we need the following special
case of \cite[III.12.11b]{hart}: 
\begin{prop}\label{killsfibers}\cite[p.52,Cor $1\frac{1}{2}$]{mumford}
Let $\rho :X\rightarrow Y$ be a flat morphism of projective
varieties, $\F$ a locally free sheaf on $X$.  If $R^i\rho_*\F=0$ for
all $i\geq i_0$, then $H^i(X_y, \F_y)=0$ for all $y\in Y$ and all
$i\geq i_0$. 
{\nopagebreak \hfill $\Box$ \par \medskip}
\end{prop}

Aside from Lemma~\ref{reduction}, the proof of Theorem~\ref{rathext} is
fairly straightforward.  We hope to clarify the idea by giving the proof now,
referencing the necessary Lemmas below:  

\begin{proof}(of Theorem~\ref{rathext})
The case $i>2$ is automatic by Lemma~\ref{elem}.  The case $i=2$
is contained in Proposition~\ref{basic}.
For $i=1$, by Proposition~\ref{basic} we need only prove the result
for $k=3$.

First note that as $L(2p-2q)$ is $4$-very ample, we can apply
Theorem~\ref{terracini} to $L(2p-2q)$ to obtain $f:X_2\rightarrow
X_0$.   Furthermore, as the restriction of $\O_{X_1}(3H-2E_1)$ to a
fiber of the $\P^1$ bundle $\wt{\Sigma_1}\rightarrow S^2C$ is
$\O_{\P^1}(-1)$, we see immediately 
\begin{eqnarray}
H^i(X_2,\O(3H-2E_1-E_2)) & = & H^i(X_1,\O(3H-2E_1))
\label{equalgroups}  \\
 & = & H^i(\P\Gamma(C,L(2p-2q)),\I_C^2(3)) \nonumber
\end{eqnarray}

Now, beginning anew with $L(2p)$,  apply Theorem~\ref{terracini} to
$L(2p)$.  This yields
$$f:\wt{X}=X_3\rightarrow X_2\rightarrow X_1\rightarrow X_0=\P
H^0(C,L(2p))$$  where $X_{i+1}=\blow{X_i}{\wt{\Sigma_i}}$.
We deduce the desired vanishing from the sequence

\begin{tabular}{rcl}
$0\rightarrow$ & $\O_{\wt{X}}(kH-4E_1-2E_2-E_3)\rightarrow$
& $\O_{\wt{X}}(kH-3E_1-2E_2-E_3)$   \\
  &   $\rightarrow \O_{E_1}(kH-3E_1-2E_2-E_3)$ & $\rightarrow 0$ 
\end{tabular}
where $k\in\Z$ is arbitrary.
The fact that $H^2(\P H^0(C,L(2p-2q)),\I_C^2(3))=0$ (this is true by
equation (\ref{equalgroups}) and Proposition~\ref{basic})
implies that $$R^2f_*\O_{E_1}(kH-3E_1-2E_2-E_3)=0$$  By
Proposition~\ref{killsfibers}, if $R^1f_*\O_{E_1}(kH-3E_1-2E_2-E_3)=0$
then the cohomology along the fibers vanishes, implying the groups in
(\ref{equalgroups}) vanish (note the higher direct images vanish by
Lemma~\ref{elem}). 

$R^1f_*\O_{\wt{X}}(kH-3E_1-2E_2-E_3)=0$ is shown in Lemma~\ref{killr1}. 

$R^2f_*\O_{\wt{X}}(kH-4E_1-2E_2-E_3)=0$ is more difficult and is shown
in Lemma~\ref{killr2}. 
\end{proof}

\begin{lemma}\label{killr1}
Under the hypotheses of Theorem~\ref{rathext}, apply
Theorem~\ref{terracini} to obtain $f:X_3\rightarrow X_0$.  Then 
$R^1f_*\O_{X_3}(kH-3E_1-2E_2-E_3)=0$.
\end{lemma}

\begin{proof}
From Lemma~\ref{tech} parts $3$ and $4$, we have $R^1f_*\O_{X_3}(kH-E_3)=0$.
Using part $3$ of Lemma~\ref{elem} to check the vanishing of $R^0f_*$ of
the rightmost term of sequences of the form 
$$\ses{\O_{X_3}(kH-E_2-E_3)}{\O_{X_3}(kH-E_3)}{\O_{E_2}(kH-E_3)}$$
and

\begin{tabular}{rcl}
$0\rightarrow$ & $\O_{X_3}(kH-E_1-2E_2-E_3)\rightarrow$ &
$\O_{X_3}(kH-2E_2-E_3)$ \\
 & $\rightarrow\O_{E_1}(kH-2E_2-E_3)$ & $\rightarrow 0$
\end{tabular}

gives $R^1f_*\O_{X_3}(kH-3E_1-2E_2-E_3)=0$.  
\end{proof}

\begin{lemma}\label{reduction}
Let $C\subset \P^n$ be a smooth curve embedded by a non-special line
bundle $L$ satisfying $(N_2)$.  Apply Theorem~\ref{terracini} to 
obtain $f:X_2\rightarrow X_0$.  Then $H^2(X_2,\O(kH-2E_1-E_2))=0$ for
$k\geq 3$. 
\end{lemma}

\begin{rem}{Remark}{notthateasy}
It may appear as if we have already shown this in the proof
of Theorem~\ref{rathext}.  Indeed, we know the result holds for
$k=3$.  However, the special circumstances involved in the
demonstration that $H^i(X_2,\O(3H-2E_1-E_2))=H^i(\I^2_C(3))$ do not
apply when $k\geq 4$.  Specifically, we need not have
$H^0(\I^2_C(k))\subseteq H^0(\I_{\Sigma_1}(k))$ for $k\geq 4$.
\end{rem}

\begin{proof}(of Lemma~\ref{reduction})
As before, equation (\ref{equalgroups}) and Proposition~\ref{basic}
imply the result for $k=3$.
Hence we show $H^2(X_2,\O(kH-E_1-E_2))=0$ for $k\geq 4$.  By restricting to 
$E_1$ and computing direct images (recall $E_1\rightarrow C$ is flat), this
immediately implies $H^2(X_2,\O(kH-2E_1-E_2))=0$ for $k\geq 4$.

Because $H^2(X_2,\O(kH-E_1))=0$, it suffices to show
$H^1(\wt{\Sigma_1},\O(kH-E_1))=0$. 
We prove $H^i(\wt{\Sigma_1},\O((4-i)H-E_1))=0$ and the result follows by a 
regularity argument (note that $\wt{\Sigma_1}$ is smooth and $\O(H)$
globally generated).   
As before, we have $H^3(\wt{\Sigma_1},\O(H-E_1))=0$ because the
restriction of $\O(H-E_1)$ to a fiber of the $\P^1$-bundle is
$\O_{\P^1}(-1)$. 
The fact that $H^1(\wt{\Sigma_1},\O(3H-E_1))=0$ follows  
immediately from projective normality and the first paragraph.

The final step is to note $\wt{\Sigma_1}\cap E_1=C\times C$ (this
follows from \cite{bertramhard}), hence we have the exact sequence 
$$\ses{\O_{\wt{\Sigma_1}}(2H-2E_1)}{\O_{\wt{\Sigma_1}}(2H-E_1)}{\O_{C\times
C}(2H-E_1)}$$ 
which we push down $f:\wt{\Sigma_1}\rightarrow S^2C$.  As
$\O_{\wt{\Sigma_1}}(2H-E_1)$ is 
trivial along the fibers, it is the pull back of a line bundle $\mathscr{L}$ on
$S^2C$.  As the restriction of $f$ to $C\times C$ is flat of degree two, 
$f_*\O_{C\times C}(2H-E_1)\cong \mathscr{L}\tensor(\O_{S^2C}\oplus
\E)$ for some  line bundle $\E$.  Therefore, $H^2(C\times
C,\O(2H-E_1))=0$ implies  $H^2(\wt{\Sigma_1},\O(2H-E_1))=0$.  It is,
however, not difficult to verify that
$\O_{C\times C}(2H-E_1)\cong \O_{C\times C}(L\bt L\tensor \I^2_{\Delta})$,
and the vanishing follows from the fact that $L$ is non-special and
very ample.
\end{proof}

As we will make three Formal Function calculations of essentially the
same type, we state an elementary result:
\begin{prop}\label{killRi}
Let $\rho :X\rightarrow Y$ be a morphism of projective
varieties; $X$ smooth, $Y$ normal.  Let $\F$ be a locally free sheaf on
$X$ and assume $y\in Y$ is a point such that:
\begin{enumerate}
\item $X_y$ is smooth 
\item $H^i(X_y,\F_y)=0$
\item $H^i(X_y,N_{X_y/X}^{*^{\tensor a}}\tensor\F_y)=0$ for all $a\geq 1$ 
\end{enumerate}
Then $R^i\rho_*\F$ is not supported at $y$.
\end{prop}

\begin{proof}
This follows by induction after tensoring the sequence:
$$\ses{S^aN^*}{\O_X/\I^{a+1}_{X_y}}{\O_X/\I^{a}_{X_y}}$$
by $\F$.
As we work over $\C$, hypothesis $3$ implies $H^i$ of the left term
vanishes.  By $2$, the Theorem on Formal Functions
\cite[III.11.1]{hart} implies the completion $(R^i\rho_*\F)^{^{\wedge}}_y=0$.
As $R^i\rho_*\F$ is coherent, the result follows (e.g. by
\cite[Ex.10.3]{atmac}). 
\end{proof}

\begin{lemma}\label{killr2}
With notation and hypotheses as in Lemma~\ref{killr1},
we have $R^2f_*\O_{X_3}(kH-rE_1-2E_2-E_3)=0$ for $r\geq 3$.
\end{lemma}

\begin{proof}
We proceed via Proposition~\ref{killRi}.  As
$R^2f_*\F := R^2f_*\O_{X_3}(kH-rE_1-2E_2-E_3)$ is supported
on $\Sigma_2\subset X_0$, we need to check three classes of fibers.

First, let $x\in \Sigma_2\setminus \Sigma_1$.  Then 
$B_x=f^{-1}(x)\cong \P H^0(C,L(2p-2Z))=\P^{n-4}$ where $Z\in S^3C$
determines the unique $3$-secant $\P^2$ containing $x$.  The
restriction of $\F$ to such a fiber is simply $\O_{\P^{n-4}}(1)$,
hence $H^2(\O_{B_x}(\F))=0$.  The conormal sequence for $B_x\subset
E_3\subset X_3$ is 
$$\ses{\O_{B_x}(-E_3)}{N^*_{B_x/X_3}}{\oplus\O_{B_x}}$$
The required vanishings follow after twisting by
$(N^*_{B_x})^{\tensor a}(\F)$. 

Let $x\in \Sigma_1\setminus C$.  Then $B_x=f^{-1}(x)\cong
\blow{\P^{n-2}}{C}$ with the embedding $C\hookrightarrow \P
H^0(C,L(2p-2Z))=\P^{n-2}$ where $Z\in S^2C$ determines the unique
secant line containing $x$.  The restriction of $\F$ to such a fiber
is $\O_{B_x}(2H-E)$ where $\picard(B_x)=\Z H+\Z E$. 
Therefore, $H^2(\O_{B_x}(\F))=0$ by projective normality of the
above embedding.  The conormal sequence for $B_x\subset
E_2\subset X_3$ is 
$$\ses{\O_{B_x}(-E_2)}{N^*_{B_x/X_3}}{\oplus\O_{B_x}}$$
and as above the required vanishing follows after twisting by
$(N^*_{B_x})^{\tensor a}(\F)$.  

If $x\in C$, then $B_x=f^{-1}(x)\cong \blow{\blow{\P^{n}}{C}}{\wt{\Sigma_1}}$ 
where $\P^{n}=\P H^0(C,L(2p-2x))$.  The restriction of $\F$ to such a
fiber is $\O_{B_x}(rH-2E_1-E_2)$ where $\picard(B_x)=\Z H+\Z E_1+\Z
E_2$.  By Lemma~\ref{reduction}, $H^2(\O_{B_x}(\F))=0$ for $r\geq 3$
and the vanishing of tensor powers of the conormal bundle follows
exactly as above. 
\end{proof}

We immediately recover Rathmann's result:
\begin{cor}\label{simplevanishing}
Assume $deg(L)\geq 2g+3$.  Then
\begin{enumerate}
\item $H^1(\P^n,\I_C^2(k))=0$ for $k\geq 3$ and $\I^2_C$ is $5$-regular
\item $\I^2_C$ is $4$-regular if and only if the Gauss-Wahl map
$\Phi_L:\wedge^2 \Gamma(L)\rightarrow \Gamma (K\otimes L^2)$ is
surjective (e.g. when $deg(L)\geq 3g+2$, \cite{belwahl})
\end{enumerate}
\end{cor}

\begin{proof}
For the first, we need only note that a line bundle of degree at least
$2g+3$ satisfies condition $(N_2)$ \cite{mgreen}.

For the second statement, we need only add that under our hypotheses,
surjectivity of $\Phi_L$ is equivalent to the vanishing
$H^2(\P^n,\I^2_C(2))=0$ \cite[1.7.3]{wahl}. 
\end{proof}

\section{Extended Vanishing}

In this section we extend Corollary~\ref{simplevanishing} to a result
suggested by Bertram \cite[4.3]{bertramvanishing}:

\begin{thm}\label{generalvanishing}
Let $C\subset \P^n$ be a smooth curve satisfying $(N_1)$ and assume
$H^1(\P^n,\I_C^2(k))=0$ for $k\geq 3$.  Then $H^i(\P^n,\I_C^a(k))=0$ for $k\geq
2a-1$, $i\geq 1$.
\end{thm}

\begin{proof}
For $i>2$ the result is again automatic by Lemma~\ref{elem}, and $i=2$
is in Proposition~\ref{basic}.

Recall $\wtp=\blow{\P^n}{C}$.  To settle the case $i=1$, let $V=\Gamma(\wtp,\O(2H-E))$ and let
$\varphi:\wtp\rightarrow\P^s$ be the morphism induced by
$\O_{\wtp}(2H-E)$.  We have the diagram:  
\begin{center}
{\begin{minipage}{1.5in}
\diagram
 & 0 \dto & 0 \dto & 0 \dto & \\
0 \rto & \varphi^*(\Omega^1_{\P^s}(1))\tensor \O_{\wtp}(H-E) \dto \rto
& V\tensor \O_{\wtp}(H-E)  \dto \rto & \O_{\wtp}(3H-2E) \dto \rto & 0 \\
0 \rto & \varphi^*(\Omega^1_{\P^s}(1))\tensor \O_{\wtp}(H) \dto \rto &
V\tensor \O_{\wtp}(H) \dto \rto & \O_{\wtp}(3H-E)  \dto \rto & 0 \\
0 \rto & \varphi^*(\Omega^1_{\P^s}(1))\tensor \O_{E}(H) \dto \rto
& V\tensor \O_{E}(H)\dto \rto & \O_{E}(3H-E) \dto \rto & 0 \\
 & 0  & 0  & 0  &
\enddiagram
\end{minipage}}
\end{center}

Twisting the entire diagram by $\O_{\wtp}(2H-E)$, we see from the top
row that showing $H^1(\I^3_C(5))=0$ is equivalent to showing
$H^1(E,\varphi^*(\Omega^1_{\P^s}(1))\tensor \O_{E}(3H-E))=0$.   
However, from the pictured diagram, it is easy to see that
$H^1(\I^2_C(3))=0$ is equivalent to the vanishing
$H^1(E,\varphi^*(\Omega^1_{\P^s}(1))\tensor \O_{E}(H))=0$. 
Therefore, twisting the last row by $\varphi^*(\Omega^1_{\P^s}(1))$
we need $H^2(E,\varphi^*(\Omega^1_{\P^s}(1))^{\tensor 2}\tensor
\O_{E}(H))=0$. 

Clearly, it suffices to prove that
the higher direct images of the blow down to the curve vanish.  As
$E\rightarrow C$ is flat, we only need vanishing along the fibers, which
are isomorphic to $\P^{n-2}$.  However, $\varphi$ maps a fiber of
$E\rightarrow C$ isomorphically to a {\em linearly embedded} subspace
$\P^{n-2}\subset \P^s$, and the vanishing follows easily.

Repeating this argument after tensoring by $\O_{\wt{P^n}}(m(2H-E))$
yields the stated result. 
\end{proof}

Analogous to Corollary~\ref{simplevanishing}, we have:
\begin{cor}\label{lesssimplevanishing}
Assume $deg(L)\geq 2g+3$.  Then
\begin{enumerate}
\item $H^1(\P^n,\I_C^a(k))=0$ for $k\geq 2a-1$
and $\I^a_C$ is $(2a+1)$-regular
\item $\I^a_C$ is $2a$-regular if and only if $\Phi_L$ is surjective
\end{enumerate}
\end{cor}

\begin{proof}
Part 1 is Theorem~\ref{rathext} applied to
Theorem~\ref{generalvanishing}.  Part 2 is the earlier statement that
$\Phi_L$ is surjective exactly when $H^2(\P^n,\I^2_C(2))=0$ applied to
Theorem~\ref{generalvanishing}. 
\end{proof}

We further have the immediate result on points:
\begin{cor}\label{points}
Let $\Gamma =C\cap H$ be a hyperplane section of a linearly normal
smooth curve of degree at least $2g+3$.  Then
$H^1(\I_{\Gamma}^a(k))=0$ for $k\geq 2a$, and the vanishing holds for
$k=2a-1$ if and only if $\Phi_L$ is surjective.
{\nopagebreak \hfill $\Box$ \par \medskip}
\end{cor} 
Note that the vanishing for $k\geq 2a$ is, more generally, true for
any set of $2n-1$ points in $\P^{n-1}$ in linearly general position by
\cite{ctv}, \cite{chandler}.

Combining \cite[Thm 1]{belwahl} with \cite[3.7]{chandler} yields the
amusing:
\begin{prop}
Let $C$ be a smooth curve embedded by a line bundle of degree $3g+1$.
Then $C$ is hyperelliptic if and only if the general hyperplane
section lies on a rational normal curve (in the hyperplane).
{\nopagebreak \hfill $\Box$ \par \medskip}
\end{prop}

The procedure detailed in Section 3 should be extendible via
Theorem~\ref{terracini} to give further vanishing statements for
higher degree embeddings.  In the very interesting cases of canonical
embeddings and higher dimensional varieties it seems that some sort of
converse (``ascending degree'') procedure must be worked out.  The
main difficulty in the canonical case is that canonical curves cannot
arise in the fibers of the blow up.  For varieties of higher
dimension, similar problems occur in that the fibers in the blow up
are copies of the original variety blown up at a point (though the
technique should at least reveal information in these cases).  A
somewhat greater obstacle is the lack of a structure theorem as strong
as Theorem~\ref{terracini}, though parts of this have been worked out
in \cite{vermeireflip1} and \cite{vermeireflip2}.

\end{article}
\end{document}